\documentclass[12pt,reqno]{amsart}
\usepackage{amscd,amssymb,amsthm,euscript}

\numberwithin{equation}{section}

\newcommand{\x}{\times}
\newcommand{\ssm}{\smallsetminus}

\newcommand{\cp}{\rtimes}

\newcommand{\iso}{\cong}

\newcommand{\cont}{\subseteq}

\newcommand{\D}{{\raise 0.09em \hbox{/}} \kern -.58em {\partial}}
\newcommand{\conv}{\lower 0.001em \hbox{*}}

\newcommand{\te}{\otimes}

\newcommand{\A}{\mathcal{A}}

\newcommand{\E}{\mathcal{E}}

\newcommand{\C}{\mathbb{C}}
\newcommand{\R}{\mathbb{R}}

\newcommand{\LL}{\mathcal{L}}

\newcommand{\T}{\mathcal{T}}

\newcommand{\F}{\mathcal{F}}
\newcommand{\Z}{\mathbb{Z}}

\newtheorem{theorem}{Theorem}[section]
\newtheorem{conjecture}[theorem]{Conjecture}

\newtheorem{proposition}[theorem]{Proposition}
\newtheorem{prop}[theorem]{Proposition}

\begin{document}
\title{A proof of the Gap Labeling Conjecture}
\author{Jerome Kaminker}
\address{Department of Mathematical Sciences, IUPUI, Indianapolis,
IN 46202-3216} 
\email{kaminker@math.iupui.edu}
\thanks{The first author was supported in part by an NSF Grant}
\thanks{The second author was supported in part by an NSERC Grant}

\author{Ian Putnam}
\address{Department of Mathematics, University of Victoria, Victoria,
BC V8W-3P4}
\email{putnam@math.uvic.ca}          

 
\begin{abstract}
  We will give a proof of the Gap Labeling Conjecture formulated by
  Bellissard, ~\cite{bellissard1}.  It makes use of a version of
  Connes' index theorem for foliations which is appropriate for
  foliated spaces, ~\cite{ moore-schochet:book}.  These  arise naturally in dynamics and are
  likely to have further applications.
\end{abstract}
\maketitle

\section{Introduction}

The ``Gap Labeling Conjecture'' as formulated by
Bellissard,~\cite{bellissard1}, is a statement about the possible gaps
in the spectrum of certain Schr\"odinger operators which arises in
solid state physics.  It has a reduction to a purely mathematical
statement about the range of the trace on a certain crossed-product
$C^{*}$-algebra, ~\cite{pedersen:book}.  By a $Cantor\ set$ we mean a compact, totally
disconnected metric space without isolated points. A group action is
$minimal$ if every orbit is dense.  
\begin{conjecture}
  Let $\Sigma$ be a Cantor set and let $\Sigma \x \Z^{n} \to \Sigma$
   be a free, minimal
   action of $\Z^{n}$ on $\Sigma$ with invariant probability measure
  $\mu$.  Let $\mu : C(\Sigma) \to \C$ and $\tau_{\mu} : C(\Sigma) \cp \Z^{n}
  \to \C$ be the traces induced by $\mu$ and denote the 
 induced maps on K-theory by the same.  Then
 one has
\begin{equation*}
  \mu(K_{0}(C(\Sigma))) = \tau_{\mu}(K_{0}(C(\Sigma) \cp \Z^{n})).
\end{equation*}
\end{conjecture}

Note that $K_{0}(C(\Sigma))$ is isomorphic to $C(\Sigma, \Z)$, the
group of integer valued continuous functions on $\Sigma$, and the
image under  $\mu$ is the subgroup of $\R$ generated by the measures of
the clopen subsets of $\Sigma$.

We will give a proof of this conjecture in the present paper.  It was
also proved independantly by J. Bellissard, R. Benedetti, and J.-M.
Gambaudo,~\cite{ bellissard:2}, and by M. Benameur and H. Oyono-Oyono,
~\cite{ oyono}.

The strategy of the proof is to use Connes' index theory for
foliations, but in the form presented in the book by Moore and
Schochet, ~\cite{ moore-schochet:book}.  In fact, this approach
underlies all three proofs, ~\cite{ bellissard:2,oyono}. Thus, one may
apply the index theorem to ``foliated spaces'' which are more general
than foliations.  These are spaces which have a cover by compatible
flow boxes as in the case of genuine foliations, except that the
transverse direction is not required to be $\R^{n}$.  In the case at
hand it is a Cantor set. 

There are two steps in the proof.  The main one is to  show that 
\begin{equation}
\label{containment}
  \tau_{\mu}(K_{0}(C(\Sigma) \cp \Z^{n})) \cont  \mu(K_{0}(C(\Sigma))).
\end{equation}
This will be carried out in Section~\ref{gap}.  The reverse
containment is easier and is proved in Section~\ref{index}.

\section{Index theory for foliated spaces}
\label{index}

We will work in a general framework based on the diagram below.
Let $\Sigma$ be a Cantor set provided with a free, minimal action of
$\Z^{n}$ and an invariant measure $\mu$, and let $X = \Sigma
\x_{\Z^{n}} \R^{n}$, be its suspension, i.e. the quotient of $\Sigma
\x \R^{n}$ by the diagonal action of $\Z^{n}$ .  There is a free
action of $\R^{n}$ on $X$ defined by $[x,w]\cdot v = [x, w+v]$.  

There is aM Morita equivalence between the $C^{*}$-algebras  associated to these group
actions, $C(\Sigma) \cp \Z^{n}$  and $C(X) \cp \R^{n}$,
~\cite{ rieffel:me}, which we will  need.  It is described in
more detail in the proof of Proposition ~\ref{prop1} below.  

Consider the diagram,
\begin{equation*}
\label{diagram}
  \begin{CD}
    K_{0}(C(\Sigma))@>{i_{*}}>>
    K_{0}(C(\Sigma) \cp \Z^{n}) @>{m.e.}>> K_{0}(C(X) \cp \R^{n})
    @<{\phi_{c}}<<
    K_{n}(C(X)) @>{ch^{(n)}}>>  \check H^{n }(X;\R)\\
    @VV{\mu}V @VV{\tau_{\mu}}V @VV{\widetilde {\tau_{\mu}}}V && @VV{ 
      C_{\mu}}V\\
    \R& =&\R& =&\R &&= &&\R
  \end{CD}
\end{equation*}

Here, the first horizontal arrow is induced by the inclusion of
$C(\Sigma)$ in $C(\Sigma) \cp \Z^{n}$, the second is provided by the
strong Morita equivalence between $C(\Sigma) \cp \Z^{n}$ and $C(X) \cp
\R^{n}$, the third is Connes' Thom isomorphism, and the fourth is the
$n^{th}$ component of the Chern character.  The first vertical arrow
is the map induced by integration against the invariant measure, the
second is the trace on $C(\Sigma) \cp \Z^{n}$ obtained from the
invariant measure on $\Sigma$, the third is induced by the trace
obtained from the associated invariant transverse measure on $X$, and
$C_{\mu}$ is the homomorphism defined via evaluation on the associated
Ruelle-Sullivan current.  We claim that this diagram commutes.  The
left squre commutes by definition of the trace, $\tau_{\mu}$.  The
other two squares will be shown to commute below.  In fact, the second
will follow by looking at the strong Morita equivalence and the third
requires  application

of the index theory of foliated spaces.

\begin{proposition}
\label{prop1}
  The  diagram
  \begin{equation}
    \begin{CD}
     K_{0}(C(\Sigma) \cp \Z^{n}) @>{m.e.}>> K_{0}(C(X) \cp \R^{n})\\ 
     @VV{\tau_{\mu}}V @VV{\widetilde {\tau_{\mu}}}V\\
     \R& =&\R
    \end{CD}
  \end{equation}
 commutes.
\end{proposition}
\begin{proof}
This is a standard fact and a proof is sketched in ~\cite{
  bellissard:3}.  We indicate a different, but related, justification here.

  The equivalence  bimodule exhibiting the strong Morita equivalence between
  $C(\Sigma) \cp \Z^{n}$ and $C(X) \cp \R^{n}$ 
   is obtained by completing
  $C_{c}(X \x \R^{n} )$,   ~\cite{ rieffel:me}.  Denote the resulting  bimodule
  by $\E$ and the  associated linking algebra, ~\cite{ brown-green-rieffel}, by
  $\A$.  Recall that $\A$
  can be viewed as being made up of  $2 \x 2$ matrices of the form
\begin{equation*}  
\begin{bmatrix}
    a & x \\
    \tilde y & b
  \end{bmatrix}
\end{equation*}
where $a \in C(\Sigma) \cp \Z^{n}$, $b \in C(X) \cp \R^{n}$, $x \in
\E$ and $\tilde y \in \E^{op}$.  This can be completed to a
 $C^{*}$-algebra, where the multiplication on the generators is given by
\begin{equation*}
  \begin{bmatrix}
    a & x \\
    \tilde y & b
  \end{bmatrix}
  \begin{bmatrix}
    a' & x' \\
    \tilde y' & b'
  \end{bmatrix}=
  \begin{bmatrix}
    aa'+ <x,\tilde y'>_{C(\Sigma) \cp \Z^{n}} & ax' + x b' \\
    \tilde y a' + b \tilde y' & bb' + <\tilde y,x'>_{C(X) \cp \R^{n}}
  \end{bmatrix}.
\end{equation*}
 The algebra $\A$ contains both $ C(\Sigma) \cp \Z^{n}$
and $ C(X) \cp \R^{n}$ as full hereditary subalgebras, hence the
inclusions, ${i_{1} : C(\Sigma)
  \cp \Z^{n} \to \A}$ and ${i_{2}:C(X) \cp \R^{n} \to \A }$, induce isomorphisms on
K-theory.  The given traces on the subalgebras give rise to a trace on $\A$
via $\tau(\begin{bmatrix}
  a & x \\
  \tilde y & b
  \end{bmatrix}) = \tau_{\mu}(a) + \tilde\tau_{\mu}(b)$.  The
  verification that this is in fact a trace requires checking that
  \begin{equation*}
\begin{split}
    \tau_{\mu}(aa'+ <x,\tilde y'>_{C(\Sigma) \cp \Z^{n}}) + \tilde
    \tau_{\mu}(bb' + <\tilde y,x'>_{C(X) \cp \R^{n}})\\ = \tau_{\mu}(a'a
    + <x',\tilde y>_{C(\Sigma) \cp \Z^{n}}) + \tilde
    \tau_{\mu}(b'b + <\tilde y',x>_{C(X) \cp \R^{n}}).
\end{split}
  \end{equation*}
This, in turn, comes down to showing that
\begin{equation*}
\begin{aligned}
  \tau_{\mu}( <x,\tilde y'>_{C(\Sigma) \cp \Z^{n}}) &= \tilde
    \tau_{\mu}( <\tilde y',x>_{C(X) \cp \R^{n}})\\
 \tau_{\mu}(<x',\tilde y>_{C(\Sigma) \cp \Z^{n}}) &= \tilde \tau_{\mu}( <\tilde y,x'>_{C(X) \cp \R^{n}}).
\end{aligned}
\end{equation*}
Each of these is a direct computation from the definitions of the
pairings and the map $\tilde\tau_{\mu}$.

It is easy
  to check that $\tau({i_{1}}_{*}(a)) = \tau_{\mu}(a)$ and
  $\tau({i_{2}}_{*}(b)) = \tilde \tau_{\mu}(b)$.  Since the
  isomorphism on K-theory induced by the strong Morita
  equivalence  is given by ${i_{2}}_{*}^{-1} {i_{1}}_{*}$, the
  result follows.
 \end{proof}

Since we will be using the theory of 
foliated spaces in the sense of Moore and Schochet, ~\cite{moore-schochet:book}, we make the
following observation about the suspension, $X$.

\begin{prop}
  The suspension, $X$, provided with its canonical $\R^{n}$-action, is
  a compact foliated space with transversal a Cantor set and invariant
  transverse measure obtained from $\mu$.
\end{prop}

We will have need of  Connes' Thom Isomorphism theorem for $C(X) \cp \R^{n}$.  It follows
from the work of 
 Fack and Skandalis, ~\cite{fack-skandalis:1981}, that the
isomorphism is induced by Kasparov product with a KK-element obtained
from the Dirac operator along the leaves of the foliated space, $X$.
Denoting the Connes Thom isomorphism by $\phi_{c} : K_{0}(C(X) \cp \R^{n}) \to
K_{n}(C(X))$, one has the following description.  
\begin{prop}
  The map $\phi_{c}$ is given by Kasparov product with the
  element 
  \begin{equation*}
    [\D ] \in KK^{n}(C(X),C(X)\cp \R^{n})
  \end{equation*}
obtained from the Dirac operator along the leaves of the foliated
space.  Thus, for an element $[E] \in K_{0}(C(X))$, one has 
\begin{equation*}
  \phi_{c}([E]) = Index^{an}([\D \te E]) \in K_{n}(C(X)\cp \R^{n})
\end{equation*}
\end{prop}
\begin{proof}
  This follows from ~\cite{fack-skandalis:1981}
\end{proof}

Finally, we are going to use the version of Connes' Foliation Index
Theorem as presented by Moore and Schochet in ~\cite{
  moore-schochet:book}.  The theorem provides a topological formula
for the result of pairing the analytic index of a leafwise elliptic
operator with the trace associated to a holonomy invariant transverse
measure.  The topological side is obtained by pairing a tangential
cohomology class with the Ruelle-Sullivan current associated to the
invariant transverse measure.  

The Ruelle-Sullivan current may be viewed as a homomorphism
\begin{equation*}
  C_{\mu} : H_{\tau}^{*}(X) \to \R,
\end{equation*}
where $H_{\tau}^{*}(X)$ is tangential
cohomology, ~\cite{ moore-schochet:book}.  This is essentially de Rham cohomology constructed from
forms that are smooth  in the leaf direction but only continuous transversally.
It is related to the \v Cech cohomology of $X$ by a natural map $r: \check H^{*}(X) \to H_{\tau}^{*}(X) $,
which in general is neither injective or surjective.  However, it
allows one to extend $C_{\mu}$ to $\check H(X)$ as $C_{\mu} \circ r$.
Moreover, for a foliated space such as $X$, there is a tangential
Chern character, $ch_{\tau}: K_{*}(C(X)) \to H_{\tau}^{*}(X)$ obtained
by applying  Chern-Weil to a leafwise connection.  It is
related to the usual Chern character via $r\circ ch = ch_{\tau}$.
With this notation at our disposal we have the following result.  

\begin{prop}
\label{indextheorem}
Let $C_{\mu}$ be the Ruelle-Sullivan current associated to the
  invariant transverse measure $\mu$ and let $ch^{(n)}$ denote the
  component of the Chern character in $\check H^{n}(X)$. Then one has
  \begin{equation*}
    \widetilde {\tau_{\mu}}(Index^{an}([\D \te E])) = C_{\mu}\circ
    r\circ ch^{(n)}([E]) 
  \end{equation*}
\end{prop}
\begin{proof}
  This is an application of the Foliation Index Theorem, ~\cite{
  connes:3, moore-schochet:book}.  By that
  theorem it is sufficient to show that the right side is what one
  obtains by pairing the index cohomology class with the
  Ruelle-Sullivan current.  In general, the index class is represented
  by the tangential form
  \begin{equation*}
ch(E)\wedge ch(\sigma(\D)) \wedge \T d(T\F \te \C) = ch(E) \wedge \hat A(T\F).    
  \end{equation*}

Note that $\hat A(T\F)$ is a polynomial in the Pontryjagin forms which are
obtained from a connection which can be chosen to be flat along the leaves, hence is
equal to $1$.  Since $r\circ ch(E) \in H_{\tau}^{*}(X)$ and taking
into account that the
homomorphism induced by the Ruelle-Sullivan current is zero except in
degree $n$ one has 
\begin{equation*}
  C_{\mu}(ch(E)\wedge ch(\sigma(\D)) \wedge \T d(T\F \te \C)) = C_{\mu}\circ
    r\circ ch^{(n)}(E),
\end{equation*}
as required.
\end{proof}

Given the above results, the commutativity of  the main diagram
 follows easily.  Indeed, the commutativity of the right hand
 rectangle is precisely the statement in Proposition
 ~\ref{indextheorem}.  

We record the following fact observed previously.

\begin{prop}
\label{easy}
  $\mu(K_{0}(C(\Sigma))) \cont  \tau_{\mu}(K_{0}(C(\Sigma) \cp \Z^{n}))$
\end{prop}

It remains to verify the other containment,
\begin{equation}
\label{containment}
 \tau_{\mu}(K_{0}(C(\Sigma) \cp \Z^{n}))  \cont \mu(K_{0}(C(\Sigma))), 
\end{equation}
which will be done in the next section.

\section{Construction of a transfer map}
\label{bott}

In this section we will provide the tool which allows the verification
of (\ref{containment}).  To accomplish this we will use the map,
described in Connes' book, ~\cite[p. 120]{connes:book}, which
associates to a clopen set in a transversal to a foliation, a
projection in its foliation algebra,
\begin{equation*}
  \alpha : K_0(C(\Sigma)) \to K_0(C(X) \cp \R^{n}). 
\end{equation*}

The modifications necessary to apply to the foliated space in question
are routine.  It will be used to relate Bott periodicity for
$C(\Sigma)$ to Connes' Thom Isomorphism for $C(X) \cp \R^{n}$.

Consider the transversal $\Sigma \x \{(\frac{1}{2}, \frac{1}{2}, \dots,
\frac{1}{2})\} \cont X$.  Let $U$ be
a clopen set of $\Sigma$ and $\chi_{U}$ its characteristic
function. We recall the description of the associated projection in
$C(X) \cp \R^{n}$.

One defines a function
\begin{equation*}
                               e_{U} : \Sigma \x [0,1] \x \R^{n} \to \R
\end{equation*}
which will yield an element of $C(X) \cp \R^{n}$ .  To this end, let $f
: \R^{n} \to \R$ be a continuous function with support in the cube of side
$\frac{1}{4} $ centered at $(\frac{1}{2}, \frac{1}{2}, \dots,
\frac{1}{2})$ and satisfying 
\begin{equation}
  \int_{\R^{n}} f(x)^{2} dx = 1.
\end{equation}
Set 
\begin{equation}
  e_{U}(x,t,s) = \chi_{U}(x)f(t-\frac{1}{2} )f(t-\frac{1}{2} -s). 
\end{equation}
Then it is easy to check that $e_{U}$ descends to a function on $X \x
\R^{n} $ that yields an element of $C(X) \cp \R^{n} $ which
satisfies  $e_{U} = e_{U}^{2} = e_{U}^{*}$.  We then set
\begin{equation}
  \label{alpha}
  \alpha(\chi_{U}) = e_{U}. 
\end{equation}
\begin{proposition}
  The function $\alpha$ induces a homomorphism
  \begin{equation}
    \alpha : K_0(C(\Sigma)) \to K_0(C(X)\cp \R^{n} ),
  \end{equation}
for which the following diagram commutes,
\begin{equation}
  \begin{CD}
    K_0(C(\Sigma)) @>{\alpha}>> K_0(C(X)\cp \R^{n} )\\
@VV{\mu}V @VV{\tilde \tau_{\mu}}V\\
\R &=& \R.
  \end{CD}
\end{equation}
\end{proposition}
\begin{proof}
  For the relation with the traces, we note that 
  \begin{equation}
    \tilde \tau_{\mu}(e_{U}) = \int_{\R^{n} } e_{U}(x,t,0) ~d\mu(x)
    ~dt = \int_{\R^{n} }
    \chi_{U} (x) f(t)f(t) ~dt ~d\mu (x) = \mu(U).  
  \end{equation}
The fact that
    $\alpha$ provides a well-defined homomorphism is straightforward.
\end{proof}

The main property of $\alpha$ is provided by the following result.
Let $\pi : \Sigma \x \R^{n} \to X$ be the quotient map.
Let $\LL$ be the union of all hyperplanes parallel to the coordinate
axis and going through points of $\Z^{n}$ and set $A = \pi (\Sigma \x
 \LL) $.  Let $j : X \setminus A \to X$ be the
inclusion of the open set $X \setminus A$ , which will induce a homomorphism $j_{*} :
C_{0}(X \setminus A) \to C(X)$.  Note that $C_{0}(X \setminus A) \iso
C_{0}(\Sigma \x (0,1)^{n}) \iso
C_{0}(\Sigma \x \R^{n})$. One now is able to relate the map $\alpha$ to Bott
periodicity and Connes' Thom isomorphism.      
\begin{proposition}
\label{diag1}
  There is a commutative diagram,
  \begin{equation}
\label{diagram1}
    \begin{CD}
       K_0(C(\Sigma)) @>{\alpha}>> K_0(C(X) \cp \R^{n})\\
     @VV{\beta}V @AA{\phi_{c}}A\\
     K_{n}(C_{0}(X \smallsetminus A)) @>{j_{*}}>> K_{n}(C(X)),
    \end{CD}
  \end{equation}
where $\phi_{c}$ is Connes' Thom Isomorphism and $\beta$ is the Bott
periodicity map.
\end{proposition}
\begin{proof}
We will deform the action $\Phi : X \x \R^{n} \to X$ as follows.  Let
$\theta_{r}: \R^{n} \to [0,1]$ be a family of continuous functions
which are periodic with respect to translation by $\Z^{n}$, with fundamental domain $[0,1]^{n}$, and on the
fundamental domain satisfy that 
\begin{itemize}
\item[i)] $\theta_{r}(\vec v) = 1$ on $[\frac{1}{4},\frac{3}{4}]^{n}$,
\item[ii)] $\theta_{r}(\vec v)$ decreases to $r$ on $\partial [0,1]^{n}$
  for $\vec v \in [0,1]^{n}
  \ssm [\frac{1}{4},\frac{3}{4}]^{n}$,
\item[iii)] $\theta_{r}(\vec v) > 0$ if $\vec v \not \in \L$.
\end{itemize}
 
Set $\Phi^{r}([z,\vec v],\vec w) = [z, \vec v-\theta_{r}(\vec v) \vec
w]$.  Here $[z,\vec v]$
denotes a point in $\Sigma \x_{\Z^{n}} \R^{n}$. Then it is easy to
check that the family $\Phi^{r}$ has the following properties.

\begin{itemize}
\item[i)] $\Phi^{1}$ is the given translation flow on $X$,
\item[ii)] $\Phi^{r}$ is the given translation flow on $[\frac{1}{4},
  \frac{3}{4}  ]^{n} \cont X$, for all   $0 \leq r \leq 1$.
\item[iii)] $\Phi^{0}$ leaves the subset $A \cont X$ pointwise fixed,
\item[iv)] $\Phi^{0}$ on $X \ssm A$ is conjugate  to $1 \x$translation on $\Sigma \x
  \R^{n}$ .
\end{itemize}

  It will be shown that the map $\alpha$ constructed with the action
  $\Phi^{0}$ agrees with $\phi_{c} j_{*} \beta$.  This proves the result.

The family can be used to define an action  $([0,1] \x X) \x \R^{n}
\to [0,1] \x X$ via the
formula 
\begin{equation*}
  \Phi(x, r) = (r, \Phi^{r}(x)).  
\end{equation*}

Consider the following commutative diagram.

\begin{equation*}
  \begin{CD}
    K_0(C(\Sigma)) @<{e_{0}}<< K_0(C(\Sigma \x [0,1])) @>{e_{1}}>> K_0(C(\Sigma))\\
    @VV{\beta}V @VV{\beta}V @VV{\beta}V \\
    K_0(C_{0}(\Sigma \x (0,1)^{n})) @<{e_{0}}<< K_0(C_{0}(\Sigma \x (0,1)^{n} \x [0,1])) @>{e_{1}}>> K_0(C_{0}(\Sigma\x (0,1)^{n}))\\
    @VV{j_{*}}V @VV{j_{*}}V @VV{j_{*}}V \\
    K_0(C(X)) @<{e_{0}}<< K_0(C(X\x [0,1])) @>{e_{1}}>> K_0(C(X))\\
    @VV{\phi_{c,0}}V @VV{\phi_{c}}V @VV{\phi_{c,1}}V \\
    K_0(C(X)\cp_{\Phi^{0}}\R^{n} ) @<{e_{0}}<< K_0(C(X \x [0,1])\cp_{\Phi}\R^{n} ) @>{e_{1}}>> K_0(C(X)\cp_{\Phi^{1}}\R^{n} ).\\
  \end{CD}
\end{equation*}

The horizontal maps are induced by evaluation at $0$ and $1$ and are
all isomorphisms.  Moreover, except for the bottom row, the
compositions $ \epsilon_{1} \epsilon_{0}^{-1}$ are the identity
homomorphism.  The vertical maps $\phi_{c,0}$, $\phi_{c}$, and
$\phi_{c,1}$ are Connes' Thom Isomorphism for the respective actions,
and $\beta$ denotes Bott periodicity.

Now,  the composition on the left side takes an element
$[\chi_{U}]$ to the element $\alpha([\chi_{U}]) = [e_{U}]_{0}$ for the
action $\Phi^{0}$.  Further, since $e_{U}$ is supported where the
actions $\Phi^{r}$ all agree, one obtains that $\epsilon_{1}
\epsilon_{0}^{-1} ([{e_{U}}_{0}]) = [{e_{U}}_{1}]$.  But, then, by
commutativity of the diagram the result follows.
\end{proof} 

\section{The gap labeling theorem}
\label{gap}

In this section we will complete
the proof of the main theorem.  Recall that we must show the
following containment holds,

  \begin{equation*}
   \tau_{\mu}(K_{0}(C(\Sigma) \cp \Z^{n})) \cont  \mu(K_{0}(C(\Sigma))) .
\end{equation*}

As a preliminary step we will look carefully at the following diagram,

\begin{equation}
  \begin{CD}
\label{diagram2}
    K_{0}(C_{0}(X \ssm A)) @>{j^{*}}>> K_{0}(C(X))\\
    @VV{ch^{(n)}}V @VV{ch^{(n)}}V\\
    \check H^{n}(X / A) @>{j^{*}}>> \check H^{n}(X).
  \end{CD}
\end{equation}

We will make two observations.
\begin{proposition}
  The map $ch^{(n)} : K_{0}(C_{0}(X \ssm A)) \to  \check H^{n}(X / A)$ is
  an isomorphism.
\end{proposition}
\begin{proof}
  The space $(X \ssm A)^{+} \iso (\Sigma \x \R^{n})^{+}$ is the inverse
  limit of finite wedges of $n$-dimensional spheres.  This is because
  $\Sigma$, being a Cantor set, is the inverse limit of finite sets.
  Since $ch^{n}$ is an isomorphism on each of the finite wedges,
  passing to the limit yields the result.
\end{proof}

\begin{proposition}
  The map $j^{*} : \check H^{n}(X / A) \to \check H^{n}(X)$ is onto.
\end{proposition}
\begin{proof}
  The map $j^{*}$ fits into the long exact sequence of the pair $(X,A)$
  and the next term is $\check H^{n}(A)$.    Recall that the definition of
  cohomological dimension of a space $X$ is
  \begin{equation}
    \dim_{\R}(X) = \sup\{ k | \check H^{k}(X,B;\R) \neq 0 \ \text{for
    some B} \cont X \}.
  \end{equation}
  Thus, it will be sufficient to show that $\dim_{\R}(A) < n$, ~\cite{
    dranishnikov-survey}.  
  
  First, we note that $A = \bigcup_{i=1}^{n} \pi (\Sigma \x \R^{n-1}_{(i)})$,
  where $\pi : \Sigma \x \R^{n} \to X$ is the projection onto the
  quotient and $\R^{n-1}_{(i)}$ denotes the points with $i^{th}$
  coordinate zero.  Now, $\pi (\Sigma \x \R^{n-1}_{(i)})$ is the total
  space of a fiber bundle with base $T^{n-1}$ and fiber a Cantor set
  $C$.  This, in turn, is a finite union of compact sets, each
  homeomorphic to $D^{n-1} \x C$, where $D^{n-1}$ is an $n-1$ disk,  and each  has
  $\dim_{\R}(D^{n-1} \x C) = n-1$.  Thus, $\dim_{\R}(\pi (\Sigma \x
  \R^{n-1}_{(i)})) = n-1$, and hence, $\dim_{\R}(A) = n-1$ since,
  again, it is a finite union of compact sets with that property.
  (See ~\cite{ dranishnikov-survey} for the properties of
  cohomological dimension needed in the above argument.)
\end{proof}

Next we assemble a larger diagram which contains ~\ref{diagram1} and
~\ref{diagram2}.

\begin{equation}
  \begin{CD}
\label{diagram3}
    K_{0}(C(\Sigma)) @>{\alpha}>> K_{0}(C(X)\cp \R^{n}) @>{\tilde
    \tau_{\mu}}>> \R\\
    @VV{\beta}V @AA{\Phi_{c}}A @VVV\\
    K^{n}(X,A) @>{j^{*}}>> K^{n}(X) && \R\\
    @VV{ch^{(n)}}V @VV{ch^{(n)}}V @AAA\\
    \check H^{n}(X,A) @>{j^{*}}>> \check H^{n}(X) @>{C_{\mu}\circ r}>> \R.
  \end{CD}
\end{equation}

The top left hand square commutes by Proposition ~\ref{diag1} and the
bottom left one does by naturality of the Chern character.  The right
hand rectangle commutes by the results in Section ~\ref{index}.  Note
that both vertical maps on the left are isomorphisms and the bottom
$j^{*}$ is onto.  We will now use this to obtain a proof of the Gap
Labeling Conjecture.
\begin{theorem}
  Let $\Z^{n}$ act minimally  on a Cantor set,
  $\Sigma$.  Consider the diagram
  \begin{equation}
    \begin{CD}
      K_{0}(C(\Sigma)) @>>> K_{0}(C(\Sigma)\cp \Z^{n})\\
      @VV{\mu}V @VV{ \tau_{\mu}}V\\
      \R& = &\R.
    \end{CD}
  \end{equation}
Then one has
\begin{equation*}
  \mu(K_{0}(C(\Sigma))) =  \tau_{\mu}(K_{0}(C(\Sigma)\cp \Z^{n})).
\end{equation*}
\end{theorem}
\begin{proof}
  It is sufficient to show that if $\lambda = \tau_{\mu}(x)$, for some
  $x \in K_{0}(C(\Sigma)\cp \Z^{n}) $, then there exists a $y \in
  K_{0}(C(\Sigma))$ with $\tau_{\mu}(x) = \mu(y)$.  To this end, we
  let $x'$ be the element of $K_{0}(C(X) \cp \R^{n})$ such that
  $\tilde \tau_{\mu}(x') = \tau_{\mu}(x)$.  We will find a $y \in
  K_{0}(C(\Sigma))$ with $\tilde \tau_{\mu}(x') = \mu(y)$.  Because
  the bottom $j^{*}$ is onto in diagram ~\ref{diagram3}, there is a $y
  \in K_{0}(C(\Sigma))$ such that $(C_{\mu}\circ r) j^{*} ch^{(n)} \beta
  (y) = \tilde \tau_{\mu}(x')$.  But by the commutativity of the
  diagram we must also have $\tilde \tau_{\mu}(\alpha(y)) = \tilde
  \tau_{\mu}(x')$.  By the basic property of $\alpha$ this yields that
  $ \mu(y) = \tilde \tau_{\mu}(x')$ which equals $\tau_{\mu}(x)$
\end{proof}

 \section{A remark on tilings and dynamics}
\label{physics}

Bellissard's original formulation of the Gap Labeling problem was for
aperiodic tiling systems.  However, we will show that the present
setting of the problem which we have addressed in this
paper---i.e. free, minimal actions of $\Z^{n}$ on Cantor sets---is
actually general enough to encompass many such tiling systems.  

 The following
proof is based on two key ingredients: a result of Sadun and Williams
and an observation which arose during a very stimulating conversation
involving Nic Ormes, Charles Radin and the second author.  For the
terminology refer to \cite{kellendonk-putnam}.

\begin{theorem}
  Suppose that $T$ is an aperiodic tiling satisfying the finite pattern condition
  and having only finitely many tile orientations.  Suppose that
  $\Omega$ is the continuous hull associated with $T$, as described in
  \cite{kellendonk-putnam}, together with the natural action of
  $\R^{n}$.  Then there is a Cantor set $\Sigma$, with a
  minimal action of $\Z^{n}$ on it,  such that $C(\Omega) \cp
  \R^{n}$ and $C(\Sigma) \cp \Z^{n}$ are strongly
  Morita equivalent.
\end{theorem}

\begin{proof}
  The result of Sadun and Williams \cite{sadun-williams} states that
  there is a Cantor set $\Sigma$ provided  with a minimal $\Z^{n}$-action such
  that the space of the suspended action, $\Sigma \times_{\Z^{n}}
   \R^{n}$, is homeomorphic to $\Omega$.  Unfortunately,
  this homeomorphism  is not a conjugacy of the $\R^{n}$ actions.
  To get around this we will bring in the fundamental groupoids, (c.f.
  ~\cite{paterson:book}), of each of these
  spaces. It is easy to see that the homeomorphism between the spaces induces an
  isomorphism between the $C^{*}$-algebras of their fundamental
  groupoids.
  
  Consider first the fundamental groupoid of $\Omega$.  We denote it by
  $\Pi(\Omega)$. There is a map of the groupoid $\Omega
  \times \R^{n}$ into $\Pi(\Omega)$ defined by sending a pair $(T,x)$ to the
  homotopy class of the path $\alpha(t) = T + tx$,  for $t \in [0,1]$.

It follows from the structure of the space $\Omega$ that this map is
an isomorphism of topological groupoids and hence induces an
isomorphism between their $C^{*}$-algebras. An analogous argument
shows that the same result holds for $\Sigma \times_{\Z^{n}} \R^{n}$.
Thus, we have $C(\Sigma) \cp \Z^{n}$ is strong Morita equivalent to $C(\Sigma
\times_{\Z^{n}} \R^{n}) \cp \R^{n}$, which is isomorphic to $C(\Omega)
\cp \R^{n}$.
\end{proof}

\end{document}